\input amstex
\magnification=\magstep1
\documentstyle{amsppt}
\advance\voffset by 0.32 truein
\advance\hoffset by 0.16 truein
\pagewidth{6.0 truein}
\pageheight{8.7 truein}
\font\smallrm=cmr8
\font\extralargebf=cmbx18
\font\sc=cmcsc10
 \font\tenbi=cmmi14
 \font\sevenbi=cmmi10 \font\fivebi=cmmi7
 \newfam\bifam \def\bi{\fam\bifam} \textfont\bifam=\tenbi
 \scriptfont\bifam=\sevenbi \scriptscriptfont\bifam=\fivebi
 \mathchardef\variablemega="7121 \def\bigomega{{\bi\variablemega}}
 \mathchardef\variablenu="7117 
\def\smashedrightarrow{\setbox0=\hbox{$\rightarrow$}\ht0=1pt\box0}
\def\risom{\buildrel{\hskip-0.08cm\sim}\over{\smashedrightarrow}}
\document

\parindent=0pt

\baselineskip=12pt

{\extralargebf Limits of Weierstrass points in regular}
\smallpagebreak

{\extralargebf smoothings of curves with two components}

\vskip0.4cm

{\bf Eduardo Esteves}

\baselineskip=10pt

{\smallrm 
Instituto de Matem\'atica Pura e Aplicada - 
Estrada Dona Castorina 110, 22460-320 Rio de Janeiro RJ, Brazil}

{\smallrm E-mail: esteves\@impa.br}

\medpagebreak

\baselineskip=12pt

{\bf Nivaldo Medeiros}

\baselineskip=10pt

{\smallrm Universidade Federal Fluminense - Instituto de Matem\'atica - 
Rua M\'{a}rio Santos Braga, s/n, Va\-lon\-guinho, 24020-005 
Niter\'{o}i RJ, Brazil}

{\smallrm E-mail: nivaldo\@impa.br}

\baselineskip=12pt

\vskip0.5cm

{\bf Abstract.} In \cite{3} D. Eisenbud and J. Harris 
posed the following question: {\it What are the limits of Weierstrass points 
in families of curves degenerating to stable curves \underbar{\rm not} 
of compact 
type?} We answer their question for one-dimensional families of smooth 
curves degenerating to stable curves with just 
two components meeting at points in general position. In this note we treat 
only those families whose total space is regular. 
Nevertheless, we announce here our most general answer, to be presented in 
detail in \cite{5}.

\parindent=8pt

\vskip0.5cm

\noindent {\bf 1. Regularly smoothable linear systems}

\vskip0.3cm 

Let $C$ be a connected, projective, nodal curve defined 
over an algebraically closed field $k$. 
Let $C_1,\dots,C_n$ be its irreducible components.
Let $B:=\text{Spec}(k[[t]])$; let $o$ denote its special point and 
$\eta$ its generic point. A projective and flat map $\pi\:S\to B$ 
is said to be a 
{\it smoothing of $C$} if the generic fiber $S_\eta$ is smooth and 
the special fiber $S_o$ is isomorphic to $C$. In addition, if 
$S$ is regular then $\pi$ is called a {\it regular smoothing}. 
If $\pi\:S\to B$ is a regular smoothing, 
then $C_1,\dots,C_n$ are Cartier divisors on $S$, and 
$C_1+\cdots+C_n\equiv 0$.

Assume from now on that $n=2$. Let $\Delta$ be the reduced Weil divisor 
with support $|\Delta|=C_1\cap C_2$ and 
$\delta:=\deg\Delta$. For $i=1,2$ let $g_i$ denote 
the arithmetic genus and $\bigomega_i$ the dualizing sheaf of $C_i$. 
Then $g:=g_1+g_2+\delta-1$ is the arithmetic genus of $C$. Let 
$\bigomega$ denote the dualizing sheaf of $C$. 

To avoid exceptional cases, assume from now on 
that $C$ is semi-stable, that is, 
assume that $\delta>1$ or 
$g_1g_2>0$. Let $\ell_i:=\lceil g_{3-i}/\delta\rceil$ 
and $m_i:=\ell_i\delta-g_{3-i}$ for $i=1,2$. 
If $\ell_1\ell_2\neq 0$ set $\lambda_i:=\ell_i/\gcd(\ell_1,\ell_2)$ for 
$i=1,2$. Let
$L_{i,j}:=\bigomega_j((1+(-1)^{i-j}\ell_i)\Delta)$ for all $i,j\in\{1,2\}$.

If $\pi$ is a regular smoothing of $C$, let $\bigomega_\pi$ be its 
(relative) dualizing sheaf. Put 
$\Cal L_{\pi,i}:=\bigomega_\pi(-\ell_iC_i)$ and 
$L_{\pi,i}:=\left.\Cal L_{\pi,i}\right|_C$ for $i=1,2$. 
Note that $L_{\pi,i}|_{C_j}\cong L_{i,j}$ for all $i,j\in\{1,2\}$. 
Once we fix isomorphisms we obtain restriction maps 
$$
\rho_{\pi,i,j}\:H^0(L_{\pi,i})\to H^0(L_{i,j})
$$
for all $i,j\in\{1,2\}$. 
Let $V_{\pi,i}:=\text{Im}(\rho_{\pi,i,i})$ for $i=1,2$. 

\vskip0.4cm

\noindent{\sc Lemma 1.} {\it For $i=1,2$ consider the following condition:
$$
h^0(\bigomega_{3-i}(-n\Delta))=\max(g_{3-i}-n\delta,0)\text{ for every }
n\geq 0.\tag{1.$i$}
$$
Let $\pi$ be a regular smoothing of $C$. 
If {\rm (1.$i$)} holds, then $h^0(L_{\pi,i})=g$ and 
$\rho_{\pi,i,j}$ is injective for $j=i$ and non-zero for 
$j=3-i$.}

\vskip0.4cm

{\it Proof.} The natural exact sequence,
$$
0 @>>> L_{\pi,i}|_{C_i}(-\Delta) @>>> L_{\pi,i} @>>> L_{\pi,i}|_{C_{3-i}} 
@>>> 0,
$$
induces an exact sequence on global sections,
$$
0 @>>> H^0(L_{i,i}(-\Delta)) @>>> H^0(L_{\pi,i}) @>\rho_{\pi,i,3-i}>> 
H^0(L_{i,3-i}).\tag{2}
$$ 
Since (2) is exact, 
$$
h^0(L_{\pi,i})\leq h^0(L_{i,i}(-\Delta))+h^0(L_{i,3-i}).
$$
Now, $h^0(L_{i,i}(-\Delta))+h^0(L_{i,3-i})=g$ 
by Riemann-Roch and (1.$i$). So $h^0(L_{\pi,i})\leq g$. 
But $L_{\pi,i}=\bigomega_\pi(-\ell_iC_i)|_C$, hence 
$h^0(L_{\pi,i})\geq g$ by semi-continuity. So $h^0(L_{\pi,i})=g$. 

Now, $g>g+m_i-\delta$, hence $h^0(L_{\pi,i})>h^0(L_{i,i}(-\Delta))$ 
if $\ell_i>0$. If $\ell_i=0$, then $g_{3-i}=0$; since $C$ is 
semi-stable, $\delta>1$ and hence $g>g_i$. It follows that 
$h^0(L_{\pi,i})>h^0(L_{i,i}(-\Delta))$ as well. In any case, 
$\rho_{\pi,i,3-i}\neq0$ by the exactness of (2).

Finally, $h^0(L_{i,3-i}(-\Delta))=0$ by (1.$i$), because 
$g_{3-i}\leq\ell_i\delta$. So $\rho_{\pi,i,i}$ is 1--1. 
The proof is complete.

\vskip0.4cm

\noindent{\sc Theorem 2.}  {\it For $i=1,2$ let $L_i$ be an invertible sheaf 
on $C$.

{\rm (a)} If $\ell_i\neq 0$, then 
there is a regular smoothing $\pi$ of $C$ such that $L_{\pi,i}\cong L_i$ 
if and only if $L_i|_{C_j}\cong L_{i,j}$ for $j=1,2$.

{\rm (b)} If $\ell_1\ell_2\neq 0$, then there is a 
regular smoothing $\pi$ of $C$ such that $L_{\pi,i}\cong L_i$ for 
$i=1,2$ if and only if $L_i|_{C_j}\cong L_{i,j}$ for all $i,j\in\{1,2\}$ 
and $L_1^{\lambda_2}L_2^{\lambda_1}\cong\bigomega^{\lambda_1+\lambda_2}$.}

\vskip0.4cm

{\it Proof.}  Let $\Bbb P$ denote the Picard group of $C$ and 
$\Bbb T\subseteq\Bbb P$ the subgroup of sheaves 
whose restrictions to $C_1$ and $C_2$ are trivial. Since $C$ is nodal, 
$\Bbb T$ is a torus of dimension $\delta-1$. Let $\Bbb E\subseteq\Bbb P$ 
denote the set of invertible sheaves $N$ on $C$ such that 
$N|_{C_i}\cong\Cal O_{C_i}((-1)^i\Delta)$ for $i=1,2$. Note 
that $\Bbb T\Bbb E\subseteq\Bbb E$.

Let's prove (a). If $\pi$ is a regular smoothing of $C$ then 
$L_{\pi,i}|_{C_j}\cong L_{i,j}$ for $j=1,2$, as observed before. 
Conversely, suppose that $L_i|_{C_j}\cong L_{i,j}$ for $j=1,2$. Let 
$N\in\Bbb E$ and 
$M:=L_i\bigomega^{-1}N^{(-1)^{i+1}\ell_i}$. Then $M\in\Bbb T$. Since 
$\Bbb T$ is a torus and $\ell_i\neq 0$, there is $Q\in\Bbb T$ such that 
$M=Q^{(-1)^i\ell_i}$. Then $L_i=\bigomega(QN)^{(-1)^i\ell_i}$. Since 
$QN\in\Bbb E$, by (\cite{6}, Prop. 3.16) there is a regular smoothing 
$\pi\:S\to B$ of $C$ such that $\Cal O_S(C_1)|_C\cong QN$. 
Then $L_{\pi,i}\cong L_i$, completing the proof of (a).

Let's prove (b). If $\pi$ is a regular smoothing of $C$ then 
$L_{\pi,i}|_{C_j}\cong L_{i,j}$ for all $i,j\in\{1,2\}$ and
$$
L_{\pi,1}^{\lambda_2}L_{\pi,2}^{\lambda_1}\cong
\bigomega_\pi^{\lambda_1+\lambda_2}
(-\lambda_2\ell_1C_1-\lambda_1\ell_2C_2)|_C.
$$
Since $\lambda_2\ell_1$ and $\lambda_1\ell_2$ are equal and 
$C_1+C_2\equiv 0$, we have $L_{\pi,1}^{\lambda_2}L_{\pi,2}^{\lambda_1}\cong
\bigomega^{\lambda_1+\lambda_2}$. 

Conversely, suppose that 
$L_i|_{C_j}\cong L_{i,j}$ for all $i,j\in\{1,2\}$ 
and $L_1^{\lambda_2}L_2^{\lambda_1}\cong\bigomega^{\lambda_1+\lambda_2}$. 
Let $N\in\Bbb E$. 
For $i=1,2$ set $M_i:=L_i\bigomega^{-1}N^{(-1)^{i+1}\ell_i}$. 
Then $M_1$ and $M_2$ are in $\Bbb T$ and 
$M_1^{\lambda_2}M_2^{\lambda_1}\cong\Cal O_C$. Since 
$\gcd(\lambda_1,\lambda_2)=1$ and $\Bbb T$ is a torus, 
there is $Q\in\Bbb T$ such that 
$M_i\cong Q^{(-1)^i\ell_i}$ for $i=1,2$. 
Then $L_i\cong\bigomega(QN)^{(-1)^i\ell_i}$ for $i=1,2$. 
Since $QN\in\Bbb E$, 
by (\cite{6}, Prop. 3.16) 
there is a regular smoothing $\pi\:S\to B$ of $C$ such that 
$\Cal O_S(C_1)|_C\cong QN$. Then $L_{\pi,i}\cong L_i$ for $i=1,2$. 
The proof is complete.

\vskip0.4cm

Note that $\ell_i\leq g_{3-i}$ and so
$L_{i,i}\subseteq\bigomega_i((1+g_{3-i})\Delta)$ for $i=1,2$. Let
$\Bbb G:=\Bbb G_1\times\Bbb G_2$ where 
$$
\Bbb G_i:=\text{Grass}_g(H^0(\bigomega_i((1+g_{3-i})\Delta)))
\text{ for $i=1,2$.}
$$

Suppose (1.$i$). We say 
that $V_i\in\Bbb G_i$ 
is {\it regularly smoothable} if there is a regular smoothing 
$\pi$ of $C$ such that $V_i=V_{\pi,i}$. Let 
$\widetilde{\Bbb V}_i\subseteq\Bbb G_i$ 
be the subset of regularly smoothable subspaces.

Suppose (1.1--2). We say that $\nu=(V_1,V_2)\in\Bbb G$ 
is {\it regularly smoothable} if there is a regular 
smoothing $\pi$ of $C$ such that $V_i=V_{\pi,i}$ for $i=1,2$, and denote 
$\nu=\nu_\pi$. Let $\widetilde{\Bbb V}\subseteq\Bbb G$ 
denote the subset of regularly smoothable 
pairs. Of course, 
$\widetilde{\Bbb V}\subseteq\widetilde{\Bbb V}_1\times\widetilde{\Bbb V}_2$.

\vskip0.4cm

\noindent{\sc Theorem 3.} 
{\it For $i=1,2$ consider the following condition:
$$
h^0(\bigomega_{3-i}(-\ell_i\Delta+I))=0\text{ for every effective divisor }
I<\Delta\text{ with }\deg I=m_i.\tag{3.$i$}
$$

{\rm (a)} If {\rm (3.$i$)} holds, then 
$\widetilde{\Bbb V}_i$ is locally closed in 
$\Bbb G_i$ and isomorphic to a torus of dimension $\delta-1$, unless 
$\delta|g_{3-i}$; in the exceptional case, 
$\widetilde{\Bbb V}_i=\{H^0(L_{i,i})\}$.

{\rm (b)} If {\rm (3.1--2)} hold and $\delta{\not|}\gcd(g_1,g_2)$, 
then $\widetilde{\Bbb V}$ is closed in 
$\widetilde{\Bbb V}_1\times\widetilde{\Bbb V}_2$ 
and isomorphic to a torus of dimension $\delta-1$.}

\vskip0.4cm

{\it Proof.}  Fix isomorphisms 
$\zeta_{i,j}\:L_{i,j}|_\Delta\risom\Cal O_\Delta$, and 
let $e_{i,j}\: H^0(L_{i,j})\to H^0(\Cal O_\Delta)$ be the induced maps 
for all $i,j\in\{1,2\}$. 
For $i=1,2$ let $L_i$ be the invertible sheaf on $C$ obtained by identifying 
$L_{i,1}$ and $L_{i,2}$ along $\Delta$ by means of $\zeta_{i,1}$ and 
$\zeta_{i,2}$. If $\ell_1\ell_2\neq 0$ choose the 
$\zeta_{i,j}$ such that 
$L_1^{\lambda_2}L_2^{\lambda_1}\cong\bigomega^{\lambda_1+\lambda_2}$. 

Assume that (3.$i$) holds. Then (1.$i$) holds as well. 
It follows from Lemma 1 and Riemann-Roch that $\dim V_{\pi,i}=g$
and $\text{codim}(V_{\pi,i},H^0(L_{i,i}))=m_i$
for every regular smoothing $\pi$ of $C$. 

If $\delta|g_{3-i}$ then $m_i=0$, and so 
$\widetilde{\Bbb V}_i=\{H^0(L_{i,i})\}$. 
Suppose now that $\delta{\not|}g_{3-i}$. 
So $\ell_i>0$, 
and thus $e_{i,i}$ is surjective 
by Riemann-Roch. 
Let $G_i:=\text{Grass}_{\delta-m_i}(H^0(\Cal O_\Delta))$. 
Taking inverse images by $e_{i,i}$ gives us a closed embedding 
$\mu_i\:G_i\to\Bbb G_i$. 
Let $W_i$ be the image of $e_{i,3-i}$. 
By (1.$i$), 
$\dim W_i=\delta-m_i$, so $W_i\in G_i$. Let $T:=H^0(\Cal O^*_\Delta)$, and 
consider its natural action on $H^0(\Cal O_\Delta)$. Let 
$\Bbb O_i$ denote the orbit of $W_i$ under the induced action of $T$ on 
$G_i$. So $\Bbb O_i\subseteq G_i$ is locally closed. 
Now, by (3.$i$), all the Pl\"ucker coordinates of $W_i$ in 
$G_i$ are non-zero. In addition, 
$W_i\neq H^0(\Cal O_\Delta)$ because $\delta{\not|}g_{3-i}$. It follows that 
the orbit map $\psi_i\: T\to\Bbb O_i$ factors through an isomorphism 
$T/k^*\risom\Bbb O_i$. So $\Bbb O_i$ is 
isomorphic to a torus of dimension $\delta-1$. Since $\ell_i\neq 0$, it 
follows from Theorem 2 that $\widetilde{\Bbb V}_i=\mu_i(\Bbb O_i)$. 
Since $\mu_i$ is an embedding, $\widetilde{\Bbb V}_i$ is locally closed in 
$\Bbb G_i$ and isomorphic to a torus of dimension $\delta-1$. So (a)
is proved.

Assume now that (3.1--2) hold, and let's prove (b). If $\delta|g_{3-i}$
then $\widetilde{\Bbb V}_i=\{H^0(L_{i,i})\}$ by (a). Hence 
$\widetilde{\Bbb V}=\widetilde{\Bbb V}_1\times\widetilde{\Bbb V}_2$ and 
$\widetilde{\Bbb V}\cong\widetilde{\Bbb V}_i$, showing (b). 
Suppose now that $\delta{\not|}g_i$ for $i=1,2$. 
Then $\ell_1\ell_2\neq 0$. Consider the subgroups
$D:=\{(t_1,t_2)\in T\times T|\  t_1^{\lambda_2}=t_2^{\lambda_1}\}$ and 
$Z:=D\cap(k^*\times k^*)$. 
Since $\lambda_1$ and $\lambda_2$ are coprime, $D$ is a subtorus of 
dimension $\delta$ of $T\times T$, 
and $Z$ is a one-dimensional subtorus of $D$. 
Let $\Bbb O$ be the orbit of $(W_1,W_2)$ under the induced action of $D$ on 
the product $G_1\times G_2$. 
Then $\Bbb O\subseteq\Bbb O_1\times\Bbb O_2$, and the orbit map 
$\psi\:D\to\Bbb O$ is the restriction to $D$ of $\psi_1\times\psi_2$. Since 
$\psi_i$ factors through an isomorphism $T/k^*\risom\Bbb O_i$ for 
$i=1,2$, then $\psi$ factors through an isomorphism $D/Z\risom\Bbb O$. 
So $\Bbb O$ is closed in $\Bbb O_1\times\Bbb O_2$ and isomorphic 
to a torus of dimension $\delta-1$. Since 
$\ell_1\ell_2\neq 0$ and
$L_1^{\lambda_2}L_2^{\lambda_1}\cong\bigomega^{\lambda_1+\lambda_2}$, 
it follows from Theorem 2 that 
$\widetilde{\Bbb V}=(\mu_1\times\mu_2)(\Bbb O)$. So $\widetilde{\Bbb V}$ is 
closed in $\widetilde{\Bbb V}_1\times\widetilde{\Bbb V}_2$ and isomorphic 
to a torus of dimension $\delta-1$. 
The proof is complete.

\vskip0.4cm

If $\pi\: S\to B$ is a smoothing of $C$, let 
$W_\pi:=\overline W_\eta\cap C$, where $W_\eta\subseteq S$ is the 
Weierstrass subscheme of the generic fiber $S_\eta$ of $\pi$. Call 
the associated Weil divisor $[W_\pi]$ a {\it limit Weierstrass divisor}. 
For each pair $\nu=(V_1,V_2)\in\Bbb G$, let 
$$
W_\nu:=W_{\nu,1}+W_{\nu,2}+g(\delta-2)\Delta,
$$
where $W_{\nu,i}$ is the ramification divisor of the linear system 
$(V_i,\bigomega_i((1+g_{3-i})\Delta))$ for $i=1,2$.

\vskip0.4cm

\noindent{\sc Theorem 4.}  {\it If $k\supseteq\Bbb Q$ and 
{\rm (1.1--2)} hold, then $[W_\pi]=W_{\nu_\pi}$ for each regular smoothing 
$\pi$ of $C$.}

\vskip0.4cm

{\it Proof.}  For $i=1,2$, since (1.$i$) holds, Lemma 1 says that 
$\Cal L_{\pi,i}$ is the extension associated to $C_i$ 
of the canonical sheaf on the generic fiber of $\pi$ ({\it see} \cite{4}, 
p. 26). By (\cite{4}, Thm. 7),
$$
[W_\pi]=\overline W_{\pi,1}+\overline W_{\pi,2}+g(g-1-\ell_1-\ell_2)\Delta,
\tag{4}
$$
where $\overline W_{\pi,i}$ is the ramification divisor 
of the linear system $(V_{\pi,i},L_{i,i})$ for $i=1,2$. Now, 
$L_{i,i}=\bigomega_i((1+\ell_i)\Delta)$, hence 
$\overline W_{\pi,i}=W_{\pi,i}-g(g_{3-i}-\ell_i)\Delta$ for $i=1,2$. 
Substituting in (4) we obtain $[W_\pi]=W_{\nu_\pi}$.

\vskip0.4cm

\noindent{\sc Corollary 5.}  {\it Assume that $k\supseteq\Bbb Q$ and 
$\delta|\gcd(g_1,g_2)$. If 
$h^0(\bigomega_i(-(g_i/\delta)\Delta))=0$ for $i=1,2$, then 
$$
[W_\pi]=W_1+W_2+(g^2-\frac{g(g+1)}{\delta})\Delta,
$$
for every regular smoothing $\pi$ of $C$, where 
$W_i$ is the ramification divisor of the complete 
system $|\bigomega_i((1+g_{3-i}/\delta)\Delta)|$ for $i=1,2$.}

\vskip0.4cm

{\it Proof.}  Since $\delta|(g_1,g_2)$, we have $\ell_i=g_{3-i}/\delta$ 
for $i=1,2$. Apply (4) to finish the proof.

\vskip0.4cm

The above corollary was stated in the case $\delta=1$ 
in (\cite{3}, Cor. 4.3).

\vskip0.3cm

\noindent {\bf 2. Smoothable linear systems} 

\vskip0.3cm

We present now the main results of 
our forthcoming \cite{5}. From now on assume that the condition,
$$
h^0(\bigomega_{3-i}(-D))=0\  \  \  \  \foldedtext\foldedwidth{6cm}
{for every effective divisor $D$ on $C_{3-i}$ with 
$|D|\subseteq|\Delta|$ and $\deg D=g_{3-i}$,}\tag{5.$i$}
$$
holds for $i=1,2$. Condition (5.$i$) is stronger than (3.$i$), 
but still is of 
general type. In \cite{5} we construct a closed subvariety 
$\Bbb V\subseteq\Bbb G$ such that $\{W_\nu|\  \nu\in\Bbb V\}$ is the set of 
limit Weierstrass divisors on $C$. This construction is described below.

First of all, every smoothing of $C$ corresponds to a regular smoothing 
of a certain semi-stable model of $C$. For each $\Delta$-tuple $\mu$ of 
positive integers consider the semi-stable model $C_\mu$ of $C$ obtained 
by splitting the branches of $C$ at each $p\in\Delta$ and connecting them 
by a chain of $\mu_p-1$ smooth rational curves. Let $\Upsilon_\mu$ denote 
the collection of irreducible components of $C_\mu$. 

Let $\widetilde\pi\:S\to B$ be a regular smoothing of $C_\mu$ and 
$\bigomega_{\widetilde\pi}$ its (relative) dualizing sheaf. For each 
$\Upsilon_\mu$-tuple $\lambda$ of integers let
$$
\Cal L_{\widetilde\pi,\lambda}:=
\bigomega_{\widetilde\pi}(\sum_{E\in\Upsilon_\mu}\lambda_EE)
\quad\text{and}\quad L_{\widetilde\pi,\lambda}:=
\Cal L_{\widetilde\pi,\lambda}|_{C_\mu}.
$$
For each $E\in\Upsilon_\mu$, consider the restriction map
$\rho_{\widetilde\pi,\lambda,E}\:H^0(L_{\widetilde\pi,\lambda})\to 
H^0(L_{\widetilde\pi,\lambda}|_E).$
Let $i=1,2$. By (\cite{4}, Thm. 1), for each regular 
smoothing $\widetilde\pi$ of $C_\mu$ there is a unique 
$\Upsilon_\mu$-tuple of integers $\lambda_i$ with $\lambda_{i,C_i}=0$ such 
that $\rho_{\widetilde\pi,\lambda_i,E}$ is 
injective for $E=C_i$ and non-zero for $E\in\Upsilon_\mu-C_i$. In \cite{5} 
we prove that $\lambda_i$ depends only on $\mu$, and that 
$h^0(L_{\widetilde\pi,\lambda_i})=g$ for every regular smoothing 
$\widetilde\pi$ of 
$C_\mu$. In fact, we give a completely numerical recipe to compute 
$\lambda_i$ in terms 
of $\mu$, and show that $\lambda_{i,E}\geq 0$ for each $E\in\Upsilon_\mu$ 
and 
$\lambda_{i,E}\leq g_{3-i}$ for each $E\in\Upsilon_\mu$ intersecting $C_i$. 
So $L_{\widetilde\pi,\lambda_i}|_{C_i}\subseteq\bigomega((1+g_{3-i})\Delta)$.

A pair $\nu=(V_1,V_2)\in\Bbb G$ is said to be
{\it $\mu$-smoothable} if there is a regular 
smoothing $\widetilde\pi$ of $C_\mu$ such that 
$V_i=\text{Im}(\rho_{\widetilde\pi,\lambda_i,C_i})$ 
for $i=1,2$. We denote $\nu=\nu_{\widetilde\pi}$. 
Let $\Bbb V_\mu\subseteq\Bbb G$ be the set of 
$\mu$-smoothable pairs. In \cite{5} we describe $\Bbb V_\mu$ in a 
way similar to the way $\widetilde{\Bbb V}$ is described in 
Theorem 2. As in Theorem 3, we prove that $\Bbb V_\mu$ is locally closed 
in $\Bbb G$ and isomorphic to a torus, and compute its dimension in terms 
of $\mu$.

Let $\Bbb U\subseteq\Bbb P_{\Bbb Q}^{\Delta}$ be the open subset 
parameterizing classes of $\Delta$-tuples of positive rational numbers, 
and give it the analytic topology. 
Changing the base of a regular smoothing $\widetilde\pi$ of $C_\mu$ 
doesn't change $L_{\widetilde\pi,\mu}$. 
Hence $\Bbb V_\mu=\Bbb V_{t\mu}$ for every 
positive integer $t$. Define $\Bbb V_\mu$ for each 
$\mu\in\Bbb U$ in the natural way. 

Let
$$
\Bbb V:=\bigcup_{\mu\in\Bbb U}\Bbb V_\mu\subseteq\Bbb G.
$$
Given $\mu,\mu'\in\Bbb U$, we show in \cite{5} that 
$\Bbb V_\mu$ and $\Bbb V_{\mu'}$ either are equal or 
don't intersect. We show also that $\Bbb V$ is 
the union of finitely many $\Bbb V_\mu$, and prove the following theorem:

\vskip0.4cm

\noindent{\sc Theorem 6.} {\it Assume {\rm (5.1--2)}. Then 
$\Bbb V$ is projective, connected and of pure dimension 
$\delta-1$.}

\vskip0.4cm

To prove the theorem we show first a remarkable relation between the 
topologies of $\Bbb U$ and $\Bbb V$. By studying boundary points of 
tori orbits in Grassmannians we prove in \cite{5} that, for each 
$\mu\in\Bbb U$,
$$
\overline{\Bbb V}_\mu=\bigcup_{\mu'\in U_\mu}\Bbb V_{\mu'},
$$
where $U_\mu\subseteq\Bbb U$ is any sufficiently small (analytic) 
open neighborhood of $\mu$. 

We show also that $\Bbb V$ is irreducible if and only if 
$g_1=g_2=1$, and compute the number of irreducible 
components of $\Bbb V$ if $\delta=2$; it is $g-\gcd(g_1+1,g_2+1)$.

Finally, let $\nu=(V_1,V_2)\in\Bbb V$. By definition, there are 
a semi-stable 
model $C_\mu$ of $C$ and a regular smoothing $\widetilde\pi$ of 
$C_\mu$ such that $\nu=\nu_{\widetilde\pi}$. Let 
$\pi$ be the smoothing of $C$ induced by $\widetilde\pi$ and 
$[W_\pi]$ the limit Weierstrass divisor on $C$. In \cite{5} we 
prove that $[W_\pi]=W_\nu$ if $k\supseteq\Bbb Q$.

\vskip0.4cm

{\sc Acknowledgements.} Previous results were obtained by L. Gatto and 
M. Coppens \cite{1}, using admissible covers. We thank 
A. Bruno, M. Coppens, D. Eisenbud, L. Gatto, J. Harris and L.~Main\`o. 
The first author was partially supported by PRONEX and 
CNPq, Proc. 300004/95-8.

\enddocument